\input amstex
\documentstyle {amsppt}
\input cyracc.def
\newfam\cyrfam
\catcode `\@=11
\font@\tencyr=wncyr10
\font@\eightcyr=wncyr8
\font@\sevencyr=wncyr7
\font@\fivecyr=wncyr5
\catcode `\@=13
\addto\tenpoint{\def\cyr{\fam\cyrfam\tencyr\cyracc}}
\addto\eightpoint{\def\cyr{\fam\cyrfam\eightcyr\cyracc}}
\textfont\cyrfam=\tencyr \scriptfont\cyrfam=\sevencyr
\NoBlackBoxes

\define\de{\delta}
\define\ep{\epsilon}
\define\ga{\gamma}
\define\la{\lambda}
\define\om{\omega}
\define\si{\sigma}

\define\De{\Delta}

\define\sh{{\cyr sh}}

\define\Qs{\roman{QSym}_k}
\define\Qsz{\roman{QSym}}
\define\Qsd{\roman{QSym}^*_k}
\define\card{\operatorname{card}}
\define\conc{\operatorname{conc}}
\define\len{\operatorname{length}}

\define\phr{\operatorname{phr}}
\define\wt{\operatorname{wt}}

\define\calg{\bold C[t_1,\dots,t_n]}
\define\calgi{\bold C[[t_1,t_2\dots]]}

\define\cmult{c_{j_1,j_2,\dots,j_m}}
\define\EE{\frak E}
\define\EEr{\frak E_r}

\define\HH{\frak H^1}

\define\BBB{1}
\define\Br{2}
\define\DS{3}
\define\DKKT{4}
\define\Eh{5}
\define\F{6}
\define\GKLLRT{7}
\define\Ge{8}
\define\Go{9}
\define\Gr{10}
\define\Hz{11}
\define\Ho{12}
\define\MaR{13}
\define\Ra{14}
\define\R{15}
\define\Ro{16}
\define\Sc{17}
\define\Sw{18}
\define\TU{19}
\define\V{20}
\define\Z{21}
\topmatter
\title Quasi-shuffle Products\endtitle
\author Michael E. Hoffman \endauthor
\address Mathematics Department, U.S. Naval Academy, Annapolis, Md. 21402
\endaddress
\email meh\@nadn.navy.mil \endemail
\keywords
Hopf algebra, shuffle algebra, quasi-symmetric functions, Lyndon word,
noncommutative symmetric functions, quantum shuffle product
\endkeywords
\subjclass Primary 16W30, 16W50, 16S80; Secondary 05E05
\endsubjclass
\abstract
Given a locally finite graded set $A$ and a commutative, associative operation
on $A$ that adds degrees, we construct a commutative multiplication $*$ on
the set of noncommutative polynomials in $A$ which we call a quasi-shuffle
product; it can be viewed as a generalization of the shuffle product $\sh$.
We extend this commutative algebra structure to a Hopf algebra
$(\frak A,*,\De)$; in the case where $A$ is the set of positive integers and
the operation on $A$ is addition, this gives the Hopf algebra of
quasi-symmetric functions.  If rational coefficients are allowed, the
quasi-shuffle product is in fact no
more general than the shuffle product; we give an isomorphism exp of the
shuffle Hopf algebra $(\frak A,\sh,\De)$ onto $(\frak A,*,\De)$.  Both the
set $L$ of Lyndon words on $A$ and their images $\{\exp(w)\mid w\in L\}$ freely
generate
the algebra $(\frak A,*)$.  We also consider the graded dual of
$(\frak A,*,\De)$.  We define a deformation
$*_q$ of $*$ that coincides with $*$ when $q=1$ and is isomorphic to the
concatenation product when $q$ is not a root of unity.  Finally, we discuss
various examples, particularly the algebra of quasi-symmetric functions
(dual to the noncommutative symmetric functions) and the algebra of
Euler sums.
\endabstract
\endtopmatter
\document
\subheading {1. Introduction} Let $k$ be a subfield of $\bold C$, and let
$A$ be a locally finite graded set.  If we think of the
graded noncommutative polynomial algebra $\frak A=k\langle A\rangle$
as a vector space over $k$, we can make it commutative $k$-algebra by
giving it the shuffle multiplication $\sh$, defined inductively by
$$
a w_1\ \sh\ b w_2 = a(w_1\ \sh\ bw_2) + b(a w_1\ \sh\ w_2)
$$
for $a,b\in A$ and words $w_1,w_2$.  The commutative $k$-algebra
$(\frak A,\sh)$ is in fact a 
polynomial algebra on the Lyndon words in $\frak A$ (as
defined in \S2 below).  If we define
$$
\De(w)=\sum_{uv=w} u\otimes v ,
$$
then $(\frak A,\sh,\De)$ becomes a commutative (but not
cocommutative) Hopf algebra, usually called the shuffle Hopf algebra;
and its graded dual is the concatenation Hopf algebra (see \cite{\R},
Chapter 1).
\par
Recently another pair of dual Hopf algebras has inspired much interest.
The Hopf algebra $\bold{Sym}$ of noncommutative symmetric functions,
introduced in \cite{\GKLLRT}, has as its graded dual the Hopf algebra
of quasi-symmetric functions \cite{\Eh,\MaR}.
In a recent paper of the author \cite{\Ho}, the algebra of
quasi-symmetric functions arose via a modification of the shuffle
product, which suggested a connection between the two pairs of Hopf
algebras.  In fact, the Hopf algebra of quasi-symmetric functions
(over $k$) is known to be 
isomorphic to the shuffle Hopf algebra on a countably infinite set
of generators (with one in each positive degree).
It is the purpose of this paper to study this Hopf algebra isomorphism
in a more general setting.  (We emphasize
that we are working over a subfield $k$ of $\bold C$; if we
instead work over $\bold Z$, there is no such isomorphism--the integral
algebra of quasi-symmetric functions is a polynomial
algebra \cite{\DS,\Sc}, but the integral shuffle algebra is not \cite{\Hz}.)
\par
More explicitly, our construction is as follows.
Suppose also that for any generators
$a,b\in A$ there is another generator $[a,b]$ so that the operation
$[\cdot,\cdot]$ is commutative, associative,
and adds degrees.  If we define a ``quasi-shuffle product'' $*$ by
$$
a w_1* b w_2 = a(w_1* bw_2) + b(a w_1* w_2)+[a,b](w_1*w_2),
$$
then $(\frak A,*)$ is a commutative and associative $k$-algebra
(Theorem 2.1 below).
In fact, as we show in \S3, $(\frak A,*,\De)$ is a Hopf algebra,
which we call the quasi-shuffle Hopf algebra corresponding to $A$ and
$[\cdot,\cdot]$.  This construction gives
the quasi-symmetric functions in the case where $A$ consists of 
one element $z_i$ in each degree $i>0$, with $[z_i,z_j]=z_{i+j}$;
this and 
other examples are discussed in \S6. We give an explicit isomorphism $\exp$
from the shuffle Hopf algebra on the generating set $A$ onto any quasi-shuffle
Hopf algebra with the same generating set (Theorems 2.5 and 3.3).
This allows us to show that any quasi-shuffle algebra on $A$ is the
free polynomial algebra on Lyndon words in $\frak A$ (Theorem 2.6).
In \S4 we take (graded) duals, giving an isomorphism $\exp^*$
from the concatenation Hopf algebra to the dual of $(\frak A,*,\De)$.
\par
In \S5 we consider a $q$-deformation $*_q$ of the quasi-shuffle product,
generalizing the quantum shuffle product as defined in \cite{\DKKT} (see
also \cite{\Gr,\Ro}).  This product coincides with the quasi-shuffle
product $*$ when $q=1$, but is noncommutative when $q\ne 1$; when $q$ is not
a root of unity, we use the theorem of Varchenko \cite{\V} to prove
that the algebra $(\frak A,*_q)$ is isomorphic to the
concatenation algebra on $A$ (Theorem 5.4).  In this case, if we
declare the elements of $A$ primitive, we get a Hopf algebra
$(\frak A,*_q,\De_q)$ isomorphic to the concatenation Hopf algebra.
\par
A construction equivalent to the quasi-shuffle algebra, but (in effect)
not assuming commutativity of the operation $[\cdot,\cdot]$, was developed
independently by F. Fares \cite{\F}.  The author thanks A. Joyal for
bringing it to his attention.
\subheading {2. The algebra structure} As above we begin with the
graded noncommutative polynomial algebra $\frak A=k\langle A\rangle$
over a subfield $k\subset\bold C$, where
$A$ is a locally finite set of generators (i.e. for each positive integer
$n$ the set $A_n$ of generators in degree $n$ is finite).  We write
$\frak A_n$ for the vector space of homogeneous elements of $\frak A$ of
degree $n$.  We shall
refer to elements of $A$ as letters, and to monomials in the letters as words.
For any word $w$ we write $\ell(w)$ for
its length (the number of letters it contains) and $|w|$ for its
degree (the sum of the degrees of its factors).  The unique word of
length 0 is 1, the empty word.
\par
Now define a new multiplication $*$ on $\frak A$ by requiring that
$*$ distribute over addition, that $1*w=w*1=w$ for any word $w$, and
that, for any words $w_1,w_2$ and letters $a,b$,
$$
aw_1*bw_2 = a(w_1*bw_2) + b(aw_1*w_2) + [a,b](w_1*w_2) ,
\tag1$$
where $[\cdot,\cdot]:\bar A\times\bar A\to\bar A$ ($\bar A=A\cup\{0\}$)
is a function satisfying
\roster
\item"{S0.}"
$[a,0]=0$ for all $a\in\bar A$;
\item"{S1.}"
$[a,b]=[b,a]$ for all $a,b\in\bar A$;
\item"{S2.}"
$[[a,b],c]=[a,[b,c]]$ for all $a,b,c\in\bar A$; and
\item"{S3.}"
Either $[a,b]=0$ or $|[a,b]|=|a|+|b|$ for all $a,b\in A$.
\endroster
\proclaim{Theorem 2.1} $(\frak A,*)$ is a commutative graded
$k$-algebra.
\endproclaim
\demo{Proof} It is enough to show that the operation $*$ is commutative,
associative, and adds degrees.
For commutativity, it suffices to show $w_1*w_2=w_2*w_1$
for any words $w_1$ and $w_2$.  We proceed by induction on $\ell(w_1)+\ell(w_2)$.
Since there is nothing to prove if either $w_1$ or $w_2$ is empty, we can assume
there are letters $a,b$ so that $w_1=au$ and $w_2=bu$.  Then (1)
together with the induction hypothesis gives
$$
w_1*w_2-w_2*w_1=[a,b](u*v)-[b,a](v*u),
$$
and the right-hand side is zero by the induction hypothesis and (S1).
Similarly, for associativity it is enough to
prove $w_1*(w_2*w_3)=(w_1*w_2)*w_3$ for any words $w_1, w_2$, and $w_3$:  this
follows from induction on $\ell(w_1)+\ell(w_2)+\ell(w_3)$ using (1) and (S2).
Finally, to show $*$ adds degrees, induct on $\ell(w_1)+\ell(w_2)$ using
(1) and (S3) to prove that $|w_1*w_2|=|w_1|+|w_2|$ for any words $w_1,w_2$.
\enddemo
If $[a,b]=0$ for all $a,b\in A$, then $(\frak A,*)$ is the shuffle
algebra as usually defined (see e.g. \cite{\R}) and we write $\sh$ for the
multiplication instead of $*$.
Suppose now that the set $A$ of letters is totally ordered.
Then lexicographic ordering gives a total order
on the words:  we put $u < uv$ for any nonempty word $v$,
and $w_1 a w_2 < w_1 b w_3$ for any letters $a< b$ and words
$w_1,w_2$, and $w_3$.
We call a word $w\ne 1$ of $\frak A$ Lyndon if $w<v$ for any nontrivial
factorization $w=uv$.  The following result is proved in Chapter 6 of
\cite{\R}; it was first obtained by Radford \cite{\Ra}.
\proclaim{Theorem 2.2} The shuffle algebra $(\frak A,\sh)$ is the free
polynomial algebra on the Lyndon words.
\endproclaim
\par
We shall define an isomorphism $\exp:(\frak A,\sh)\to (\frak A,*)$.
To do so, we must first develop some notation relating to the
operation $[\cdot,\cdot]$ and compositions.  Define
inductively $[S]\in\bar A$ for any finite sequence $S$ of elements of $A$
by setting $[a]=a$ for $a\in A$, and $[a,S] = [a,[S]]$ for any $a\in A$
and sequence $S$ of elements of $A$.  
\proclaim{Proposition 2.3} (i) If $[S]=0$, then $[S']=0$ whenever
$S$ is a subsequence of $S'$;
\newline
(ii) $[S]$ does not depend on the order of
the elements of $S$ (i.e., it depends only on the underlying multiset
of $S$);
\newline
(iii) For any sequences $S_1$ and $S_2$, $[S_1\sqcup S_2]=[[S_1],[S_2]]$,
where $S_1\sqcup S_2$ denotes the concatenation of sequences $S_1$ and
$S_2$;
\newline
(iv) If $[S]\ne 0$, then the degree of $S$ is the sum of the degrees of the
elements of $S$.
\endproclaim
\demo{Proof} (i),(ii),(iii),(iv) follow from (S0),(S1),(S2),(S3) respectively.
\enddemo
\par
A composition of a positive integer $n$ is a sequence $I=(i_1,i_2,\dots, i_k)$
of positive integers such that $i_1+i_2+\dots+i_k=n$.  We call $n=|I|$ the
weight of $I$ and $k=\ell(I)$ its length; we write $\Cal C(n)$ for the set
of compositions of $n$, and $\Cal C(n,k)$ for the set of compositions of
$n$ of length $k$.  For $I\in\Cal C(n,k)$ and $J\in\Cal C(k,l)$, the
composition $J\circ I\in \Cal C(n,l)$ is given by
$$
J\circ I=(i_1+\dots+i_{j_1},i_{j_1+1}+\dots+i_{j_1+j_2},\dots,
i_{j_1+\dots+j_{l-1}+1}+\dots+i_k) .
$$
If $K=J\circ I$ for some $J$, we call $I$ a refinement of $K$ and write
$I\succcurlyeq K$.
Compositions act on words via $[\cdot,\cdot]$ as follows.  For any word
$w=a_1a_2\cdots a_n$ and composition $I=(i_1,\dots,i_l)\in\Cal C(n)$, set
$$
I[w]=[a_1,\dots,a_{i_1}][a_{i_1+1},\dots,a_{i_1+i_2}]\cdots
[a_{i_1+\dots+i_{l-1}+1},\dots,a_n] .
$$
(This is really an action in the sense that $I[J[w]]=I\circ J[w]$.)
\par
Now let $\exp:\frak A\to\frak A$ be the linear map with $\exp(1)=1$ and
$$
\exp(w)=\sum_{(i_1,\dots,i_l)\in\Cal C(\ell(w))}\frac1{i_1!\cdots i_l!}
(i_1,\dots,i_l)[w]
$$
for any nonempty word $w$ (so, e.g. $\exp(a_1a_2a_3)=a_1a_2a_3 + \frac12 [a_1,a_2]a_3 +\frac12 a_1
[a_2,a_3] + \frac16 [a_1,a_2,a_3]$).  There is an inverse $\log$ of
$\exp$ given by
$$
\log(w)=\sum_{(i_1,\dots,i_l)\in\Cal C(\ell(w))} \frac{(-1)^{\ell(w)-l}}
{i_1\cdots i_l}(i_1,\dots,i_l)[w] 
$$
for any word $w$, and extended to $\frak A$ by linearity; this follows
by taking $f(t)=e^t-1$ in the following lemma.
\proclaim{Lemma 2.4} Let $f(t)=a_1t+a_2t^2+a_3t^3+\cdots$ be a function
analytic at the origin, with $a_1\ne 0$ and $a_i\in k$ for all $i$, and
let $f^{-1}(t)=b_1t+b_2t^2+b_3t^3+\cdots$ be the inverse of $f$.
Then the map $\Psi_f:\frak A\to\frak A$ given by
$$
\Psi_f(w)=\sum_{I\in\Cal C(\ell(w))}a_{i_1}a_{i_2}\cdots a_{i_l}I[w]
$$
for words $w$, and extended linearly, has inverse $\Psi_f^{-1}=\Psi_{f^{-1}}$
given by
$$
\Psi_{f^{-1}}(w)=\sum_{I\in\Cal C(\ell(w))}b_{i_1}b_{i_2}\cdots b_{i_l}I[w] .
$$
\endproclaim
\demo{Proof} It suffices to show that $\Psi_{f^{-1}}(\Psi_f(w))=w$ for any
word $w$ of length $n\ge 1$
(Note that $\Psi_f(\Psi_{f^{-1}}(w))=w$ is then automatic, since $\Psi_f$
and $\Psi_{f^{-1}}$
can be thought of as linear maps of the vector space with basis
$\{I[w]\mid I\in\Cal C(n)\}$.)
Now for any $K=(k_1,\dots,k_l)\in\Cal C(n)$, the coefficient of $K[w]$ in
$\Psi_{f^{-1}}(\Psi_f(w))$ is
$$
\sum_{J\circ I=K}b_{j_1}b_{j_2}\cdots b_{j_l}a_{i_1}a_{i_2}\cdots a_{i_{|J|}} .
\tag2$$
We must show that (2) is 1 if $K$ is a sequence of $n$ 1's, and 0 otherwise.
To see this, let $t_1,t_2,\dots$ be commuting variables.  Then
(2) is the coefficient of $t_1^{k_1} t_2^{k_2}\cdots t_l^{k_l}$ in
$$
t_1t_2\cdots t_l=f^{-1}(f(t_1))f^{-1}(f(t_2))\cdots f^{-1}(f(t_l)) .
$$
\enddemo
\proclaim{Theorem 2.5} $\exp$ is an isomorphism of $(\frak A,\sh)$ onto
$(\frak A,*)$ (as graded $k$-algebras).
\endproclaim
\demo{Proof} From the lemma, $\exp$ is invertible.  Also, it follows from
2.3(iv) that $\exp$ preserves degree. To show $\exp$ a homomorphism it
suffices to
show $\exp(w\ \sh\ v)=\exp(w)*\exp(v)$ for any words $w,v$.  Let $w=a_1\cdots a_n$
and $v=b_1\cdots b_m$.  Evidently both $\exp(w\ \sh\ v)$ and $\exp(w)*\exp(v)$
are sums of rational multiples of terms
$$
[S_1\sqcup T_1][S_2\sqcup T_2]\cdots[S_l\sqcup T_l]
\tag3$$
where the $S_i$ and $T_i$ are subsequences of $a_1,\dots,a_n$ and
$b_1,\dots,b_m$ respectively such that
\roster
\item"{i.}"
for each $i$, at most one of $S_i$, $T_i$ is empty; and
\item"{ii.}"
the concatenation $S_1\sqcup S_2\sqcup\cdots\sqcup S_l$
is the sequence $a_1,\dots,a_n$, and similarly the $T_i$ concatenate
to give the sequence $b_1,\dots,b_m$.
\endroster
Now the term (3) arises in $\exp(w)*\exp(v)$ in only one way, and its
coefficient is
$$
\frac1{(\len S_1)!(\len S_2)!\cdots(\len S_l)!
(\len T_1)!(\len T_2)!\cdots(\len T_l)!} .
$$
On the other other hand, (3) can arise in $\exp(w\ \sh\ v)$ from
$$\multline
\binom{\len S_1\sqcup T_1}{\len S_1}
\binom{\len S_2\sqcup T_2}{\len S_2}\cdots
\binom{\len S_l\sqcup T_l}{\len S_l}=\\
\frac{(\len S_1\sqcup T_1)!\cdots(\len S_l\sqcup T_l)!}
{(\len S_1)!\cdots(\len S_l)!(\len T_1)!\cdots(\len T_l)!}
\endmultline$$
distinct terms of the shuffle product $w\ \sh\ v$, and after application
of $\exp$ each such term acquires a coefficient of
$$
\frac1{(\len S_1\sqcup T_1)!\cdots(\len S_l\sqcup T_l)!} .
$$
\enddemo
\par
It follows from Theorems 2.2 and 2.5 that $(\frak A,*)$ is the free polynomial
algebra on the elements $\{\exp(w)\mid w\ \text{is a Lyndon word}\}$.  In
fact the following is true.
\proclaim{Theorem 2.6} $(\frak A,*)$ is the free polynomial algebra on
the Lyndon words.
\endproclaim
\demo{Proof} It suffices to show that any word $w$ can be written
as a $*$-polynomial of Lyndon words.  We proceed by induction on
$\ell(w)$.  If $\ell(w)=1$ the result is immediate, since every letter is
a Lyndon word.  Now let $\ell(w)>1$: by Theorem 2.5 there are Lyndon
words $w_1,\dots, w_n$ and a polynomial $P$ so that
$$
w=P(\exp(w_1),\dots,\exp(w_n))
$$
in $(\frak A,*)$.  Note that since $\log(w)=P(w_1,\dots,w_n)$
in $(\frak A,\sh)$, we can assume every term of $P(w_1,\dots,w_n)$
(as a $\sh$-polynomial)
has length at most $\ell(w)$, since the shuffle product preserves lengths.
But then in $(\frak A,*)$,
$$
w-P(w_1,\dots,w_n)=P(\exp(w_1),\dots,\exp(w_n))-P(w_1,\dots,w_n)
$$
must consist of terms of length less than $\ell(w)$,
and so is expressible in terms of Lyndon words by the induction hypothesis.
\enddemo
\par
By the preceding result, the number of generators of $(\frak A,*)$ in
degree $n$ is the number $L_n$ of Lyndon words of degree $n$.  This
number can be calculated from Poincar\'e series
$$
A(x)=\sum_{n\ge 0}(\dim\frak A_n) x^n=\frac1{1-\sum_{n\ge 1}(\card A_n)x^n}
$$
of $\frak A$ as follows.
\proclaim{Proposition 2.7} The number $L_n$ of Lyndon words in $\frak A_n$
is given by
$$
L_n = \frac1{n}\sum_{d|n}\mu\left(\frac{n}{d}\right)c_d ,
$$
where the numbers $c_n$ are defined by
$$
x\frac{d}{dx}\log A(x)=\sum_{n\ge 1}c_n x^n 
$$
for $A(x)$ as above.
\endproclaim
\demo{Proof} In view of Theorems 2.2 and 2.6, we must have
$$
A(x)=\prod_{n\ge 1}(1-x^n)^{-L_n} .
$$
The conclusion then follows from
taking logarithms, differentiating, and using the M\"obius inversion
formula.
\enddemo
\subheading {3. The Hopf algebra structure} For basic definitions and
facts about Hopf algebras see \cite{\Sw}.  We define a comultiplication
$\De:\frak A\otimes\frak A\to\frak A$ and counit $\ep:\frak A\to k$ by
$$
\De(w)=\sum_{uv=w} u\otimes v
$$
and
$$
\ep(w)=\cases 1,&w=1\\ 0,&\text{otherwise}\endcases
$$
for any word $w$ of $\frak A$.
Then $(\frak A,\De,\ep)$ is evidently a (non-cocommutative) coalgebra.
In fact the following result holds.
\proclaim{Theorem 3.1} $\frak A$ with the $*$-multiplication and
$\De$-comultiplication is a bialgebra.
\endproclaim
\demo{Proof} It suffices to show that $\ep$ and $\De$ are $*$-homomorphisms.
The statement for $\ep$ is obvious; to show $\De(w_1)*\De(w_2)=\De(w_1*w_2)$
for any words $w_1,w_2$ use induction on $\ell(w_1)+\ell(w_2)$.
Since the result is immediate if $w_1$ or $w_2$ is
1, we can write $w_1=au$ and $w_2=bv$ for letters $a,b$ and words $u,v$.
Adopting Sweedler's sigma notation \cite{\Sw}, we write
$$
\De(u)=\sum u_{(1)}\otimes u_{(2)},\quad\text{and}\quad
\De(v)=\sum v_{(1)}\otimes v_{(2)} .
$$
Then from the definition of $\De$,
$$
\De(w_1)=\sum au_{(1)}\otimes u_{(2)}+1\otimes au
\quad\text{and}\quad
\De(w_2)=\sum bv_{(1)}\otimes v_{(2)}+1\otimes bv ,
$$
so that $\De(w_1)*\De(w_2)$ is 
$$
\sum (au_{(1)}*bv_{(1)})\otimes(u_{(2)}*v_{(2)})
+\sum au_{(1)}\otimes(u_{(2)}*bv)
+\sum bv_{(1)}\otimes(au*v_{(2)})
+1\otimes(au*bv) .
$$
Using (1), this is
$$\multline
\sum a(u_{(1)}*bv_{(1)})\otimes(u_{(2)}*v_{(2)})+
\sum b(au_{(1)}*v_{(1)})\otimes(u_{(2)}*v_{(2)})+\\
\sum [a,b](u_{(1)}*v_{(1)})\otimes(u_{(2)}*v_{(2)})+
\sum au_{(1)}\otimes(u_{(2)}*bv)+
\sum bv_{(1)}\otimes(au*v_{(2)})\\
+1\otimes a(u*w_2)+1\otimes b(w_1*v)+1\otimes[a,b](u*v),
\endmultline$$
or, applying the induction hypothesis,
$$\multline
(a\otimes 1)(\De(u)*\De(w_2))+1\otimes a(u*w_2)+
(b\otimes 1)(\De(w_1)*\De(v))+1\otimes b(w_1*v)\\
+([a,b]\otimes 1)\De(u*v)+1\otimes[a,b](u*v),
\endmultline$$
which can be recognized as $\De(w_1*w_2)=\De(a(u*w_2)+b(w_1*v)+[a,b](u*v))$.
\enddemo
Since both $*$ and $\De$ respect the grading, it follows automatically
that $\frak A$ is a Hopf algebra (cf. Lemma 2.1 of \cite{\Eh}).
In fact there are two explicit formulas for the antipode, whose agreement
is of some interest.
\proclaim{Theorem 3.2} The bialgebra $\frak A$ has antipode $S$ given
by
$$\multline
S(w) = \sum_{(i_1,\dots,i_l)\in\Cal C(n)}
(-1)^l a_1\cdots a_{i_1}*a_{i_1+1}\cdots a_{i_1+i_2}*\cdots*a_{i_1+\dots+
i_{l-1}+1}\cdots a_n \\
= (-1)^n \sum_{I\in\Cal C(n)}I[a_n a_{n-1}\cdots a_1]
\endmultline$$
for any word $w=a_1a_2\cdots a_n$ of $\frak A$.
\endproclaim
\demo{Proof} We can compute $S$ recursively from $S(1)=1$ and
$$
S(w)= - \sum_{k=0}^{n-1} S(a_1\cdots a_k)*a_{k+1}\cdots a_n
\tag4$$
for a word $w=a_1\cdots a_n$.  The first formula for $S$ then follows
easily by induction on $n$.  For the the second formula, we also proceed
by induction on $n$, following the proof of Proposition 3.4 of \cite{\Eh}.
For $w=a_1\cdots a_n$, $n>0$, the induction hypothesis and (4) give
$S(w)$ as
$$\multline
\sum_{k=0}^{n-1} \sum_{(i_1,\dots,i_l)\in\Cal C(k)}
(-1)^{k+1} (i_1,\dots,i_l)[a_k a_{k-1}\cdots a_1]*a_{k+1}\cdots a_n =\\
\sum_{k=0}^{n-1} \sum_{(i_1,\dots,i_l)\in\Cal C(k)}
(-1)^{k+1} [a_k,a_{k-1},\dots,a_{k-i_1+1}]\cdots [a_{i_l},\dots,a_1]*
a_{k+1}\cdots a_n
\endmultline$$
Now the first factor of each term of the $*$-product in the inner sum is,
from consideration of (1), one of three generators: $[a_k,\dots, a_{k-i_1+1}]$,
$[a_{k+1},a_k,\dots,a_{k-i_1+1}]$, or $a_{k+1}$.  We say the term is of
type $k$ in the first case, and of type $k+1$ in the latter two cases.
Now consider a word that appears in the expansion of $S(w)$.  If it has
type $i\le n-1$, then it occurs for both $k=i$ and $k=i-1$, and the
two occurrences will cancel.  The only words that do not cancel are those of
type $n$, which occur only for $k=n-1$: these will all carry the coefficient
$(-1)^n$, and give the second formula for $S(w)$.
\enddemo
\demo{Remark} In the case of the shuffle algebra (i.e., where $[\cdot,\cdot]$ is
identically zero), the second formula for the antipode is simply
$S(w)=(-1)^{\ell(w)}\bar w$.  Cf. \cite{\R, p. 35}.
\enddemo
\proclaim{Theorem 3.3} $\exp:\frak A\to\frak A$ is a Hopf algebra isomorphism
of $(\frak A,\sh,\De)$ onto $(\frak A,*,\De)$.
\endproclaim
\demo{Proof} We have already shown that $\exp$ is an algebra homomorphism.
It suffices to show that $\exp\circ\ep(w)=\ep\circ\exp(w)$ and
$\De\circ\exp(w)=(\exp\otimes\exp)\circ\De(w)$ for any word $w$.
The first equation is immediate, and the second follows since both sides
are equal to
$$
\sum_{uv=w}\sum\Sb (i_1,\dots,i_k)\in\Cal C(\ell(u))\\
(j_1,\dots,j_l)\in\Cal C(\ell(v))\endSb \frac1{i_1!\cdots i_k!}I[u]\otimes
\frac1{j_1!\cdots j_l!}J[v] .
$$
\enddemo
\subheading {4. Duality} The graded dual $\frak A^*=\bigoplus_{n\ge 0}\frak A_n^*$
has a basis consisting of elements $w^*$, where $w$ is a word of $\frak A$:
the pairing $(\cdot,\cdot):\frak A\otimes \frak A^*\to k$ is given by
$$
(u,v^*)=\cases 1&\text{if $u=v$}\\0&\text{otherwise.}\endcases
$$
Then the transpose of $\De$ is the concatenation
product $\conc(u^*\otimes v^*)=(uv)^*$, and the transpose of 
$\sh$ is the comultiplication $\de$ defined by
$$
\de(w^*)=\sum_{\text{words $u,v$ of $\frak A$}}(u\ \sh\ v,w^*)u^*\otimes v^* .
$$
Since $(\frak A,\sh,\De)$ is a Hopf algebra, so is its graded dual
$(\frak A^*, \conc, \de)$, which is
called the concatenation Hopf algebra in \cite{\R}.  Dualizing
$(\frak A,*,\De)$, we also have a
Hopf algebra $(\frak A^*, \conc, \de')$, where $\de'$ is the comultiplication
defined by
$$
\de'(w^*)=\sum_{\text{words $u,v$ of $\frak A$}}(u*v,w^*)u^*\otimes v^* .
$$
Then from our earlier results we have the following.
\proclaim{Theorem 4.1} There is a Hopf algebra isomorphism $\exp^*$ from
$(\frak A^*,\conc,\de')$ to $(\frak A^*,\conc,\de)$.
\endproclaim
$\exp^*$ is the transpose of $\exp$:  explicitly, $\exp^*$ is the endomorphism
of $(\frak A^*,\conc)$ with
$$
\exp^*(a^*)=\sum_{n\ge 1}\frac1{n!}\sum_{(n)[w]=a}w^*=
\sum_{n\ge 1}\sum_{[a_1,\dots,a_n]=a}\frac1{n!}(a_1\cdots a_n)^*
$$
for $a\in A$.  It has inverse $\log^*$ given by
$$
\log^*(a^*)=\sum_{n\ge 1}\frac{(-1)^{n-1}}{n}\sum_{(n)[w]=a}w^*,
\quad a\in A .
\tag5$$
The set of Lie polynomials in $\frak A^*$ is the smallest sub-vector-space
of $\frak A^*$
containing the set of generators $\{a^*\mid a\in A\}$ and closed under the
Lie bracket
$$
[P,Q]_{\roman{Lie}}=PQ-QP .
$$
Since the Lie polynomials are exactly the primitives for $\de$ \cite{\R,
Theorem 1.4}, we have the following result.
\proclaim{Theorem 4.2} The primitives for $\de'$ are elements of the
form $\log^*P$, where $P$ is a Lie polynomial.
\endproclaim
We note that $(\frak A^*,\conc,\de')$ has antipode
$$
S^*(w^*)=\sum_{v\in \Cal P(\bar w)} (-1)^{\ell(v)}v^* ,
$$
where $\bar w$ is the reverse of $w$ (i.e. $\bar w = a_n a_{n-1}\cdots a_1$
if $w=a_1 a_2\cdots a_n$) and $\Cal P(w)=\{v\mid\text{$I[v]=w$ for some $I\in
\Cal C(\ell(v))$}\}$.
\subheading {5. $q$-deformation} We now define a deformation of $(\frak A,*)$.
We again start with the noncommutative polynomial algebra $\frak A=k\langle A\rangle$
and define, for $q\in k$, a new multiplication $*_q$ by requiring that $*_q$
distribute over addition, that $w*_q1=1*_qw=w$ for any word $w$ and that
$$
aw_1*_qbw_2=a(w_1*_qbw_2)+q^{|aw_1||b|}b(aw_1*_qw_2)+q^{|w_1||b|}[a,b](w_1*_qw_2)
\tag6$$
for any words $w_1$, $w_2$ and letters $a$, $b$.
\proclaim{Theorem 5.1} $(\frak A,*_q)$ is a graded $k$-algebra, which
coincides with $(\frak A,*)$ when $q=1$.
\endproclaim
\demo{Proof} The argument is similar to that for Theorem 2.1.  It is easy
to show that $|w_1*_qw_2|=|w_1|+|w_2|$ for any words $w_1$, $w_2$ by
induction on $\ell(w_1)+\ell(w_2)$.  To show
the operation $*_q$ associative, it suffices to show that
$w_1*_q(w_2*_qw_3)=(w_1*_qw_2)*_qw_3$ for any words $w_1$, $w_2$, and $w_3$,
which we do by induction on $\ell(w_1)+\ell(w_2)+\ell(w_3)$.  We can
assume $w_i=a_iu_i$ for letters $a_i$ and words $u_i$, $i=1,2,3$.
Then $w_1*_q(w_2*_qw_3)$ is
$$\multline
a_1(u_1*_qa_2(u_2*_qw_3))+
q^{|w_1||a_2|}a_2(w_1*_q(u_2*_qw_3))+
q^{|u_1||a_2|}[a_1,a_2](u_1*_q(u_2*_qw_3))\\
+q^{|w_2||a_3|}a_1(u_1*_qa_3(w_2*_qu_3))+
q^{|w_2||a_3|+|w_1||a_3|}a_3(w_1*_q(w_2*_qu_3))+\\
q^{|w_2||a_3|+|u_1||a_3|}[a_1,a_3](u_1*_q(w_2*_qu_3))+
q^{|u_2||a_3|}a_1(u_1*_q[a_2,a_3](u_2*_qu_3))+\\
q^{|u_2||a_3|+|w_1||a_2a_3|}[a_2,a_3](w_1*_q(u_2*_qu_3))+
q^{|u_2||a_3|+|u_1||a_2a_3|}[a_1,a_2,a_3](u_1*_q(u_2*_qu_3)),
\endmultline$$
while $(w_1*_qw_2)*_qw_3$ is
$$\multline
a_1((u_1*_qw_2)*_qw_3)+
q^{|w_1w_2||a_3|}a_3(a_1(u_1*_qw_2)*_qu_3)+
q^{|u_1w_2||a_3|}[a_1,a_3]((u_1*_qw_2)*_qu_3)\\
+q^{|w_1||a_2|}a_2((w_1*_qu_2)*_qw_3)+
q^{|w_1||a_2|+|w_1w_2||a_3|}a_3(a_2(w_1*_qu_2)*_qu_3)+\\
q^{|w_1||a_2|+|w_1u_2||a_3|}[a_2,a_3]((w_1*_qu_2)*_qu_3)+
q^{|u_1||a_2|}[a_1,a_2]((u_1*_qu_2)*_qw_3)+\\
q^{|u_1||a_2|+|w_1w_2||a_3|}a_3([a_1,a_2](u_1*_qu_2)*_qu_3)+
q^{|u_1||a_2|+|u_1u_2||a_3|}[a_1,a_2,a_3]((u_1*_qu_2)*_qu_3) .
\endmultline$$
Applying the induction hypothesis, the difference is
$$\multline
a_1(u_1*_q(a_2(u_2*_qw_3)+q^{|w_2||a_3|}a_3(w_2*_qu_3)+
q^{|u_2||a_3|}[a_2,a_3](u_2*_qu_3)))\\
+q^{(|w_2|+|w_1|)|a_3|}a_3(w_1*_q(w_2*_qu_3))
-a_1((u_1*_qw_2)*_qw_3)\\
-q^{|w_1w_2||a_3|}a_3((a_1(u_1*_qw_2)
+q^{|w_1||a_2|}a_2(w_1*_qu_2)+q^{|u_1||a_2|}[a_1,a_2](u_1*_qu_2))*_qu_3),
\endmultline$$
which by application of (6) and the induction hypothesis is seen to be
zero.
\enddemo
\demo{Remark} The author arrived at the definition (6) as follows.  Knowing
the first two terms on the right-hand side from the definition of
the quantum shuffle product, he tried an arbitrary power of $q$ on the
third term, and found that the resulting product was only associative
when the exponent is as in (6).  Shortly afterward he discussed this
with J.-Y. Thibon, who directed him to \cite{\TU}, where
the rule (6) appears in the special case of the quasi-symmetric functions
(see Example 1 below).
\enddemo
\par
Of course, for $q\ne 1$ the algebra $(\frak A,*_q)$ is no longer
commutative.  For each fixed $q$, there is a homomorphism $\Phi_q$ of
graded associative $k$-algebras from the concatenation algebra
$(\frak A,\conc)$ to $(\frak A,*_q)$ defined by
$$
\Phi_q(a_1a_2\cdots a_n)=a_1*_q a_2*_q\cdots *_q a_n
$$
for letters
$a_1$, $a_2,\dots, a_n$; we call $q$ generic if $\Phi_q$ is an isomorphism
(i.e., if it is surjective).  To give an explicit formula for $\Phi_q$,
we introduce some notation.  For a permutation $\si$ of $\{1,2,\dots,n\}$,
let $\iota(\si)=\{(i,j)\mid\text{$1\le i<j\le n$ and $\si(i)>\si(j)$}\}$
be the set of inversions of $\si$, and let $C(\si)$ be the descent composition
of $\si$, i.e. the composition $(i_1,i_2,\dots,i_l)\in\Cal C(n)$ with 
$$
\si(i_1+\dots+i_{j-1}+1)<\si(i_1+\dots+i_{j-1}+2)<\cdots<
\si(i_1+\dots+i_j)
$$
for $j=1,2,\dots, l$ and $l$ minimal.  (Equivalently, $C(\si)=(i_1,\dots,i_l)$
is the composition such that the associated subset $\{i_1,i_1+i_2,\dots,
i_1+\dots+i_{l-1}\}$ of $\{1,2,\dots,n-1\}$ is the descent set of $\si$,
i.e. the set of $1\le i\le n-1$ such that $\si(i)>\si(i+1)$.)
\proclaim{Lemma 5.2} For any letters $a_1,a_2,\dots, a_n$,
$$
\Phi_q(a_1a_2\cdots a_n)=\sum_{\text{permutations $\si$}}
q^{\sum_{(i,j)\in\iota(\si)}|a_i||a_j|}
\sum_{I\succcurlyeq C(\si)}I[a_{\si(1)}a_{\si(2)}\cdots a_{\si(n)}]
$$
\endproclaim
\demo{Proof} We proceed by induction on $n$, the case $n=2$ being immediate.
Assuming the induction hypothesis, we have
$$
\Phi_q(a_1\cdots a_{n+1})=\sum_{(\si,I)\in P(n)}
q^{\sum_{(i,j)\in\iota(\si)}|a_i||a_j|}
I[a_{\si(1)}a_{\si(2)}\cdots a_{\si(n)}]*a_{n+1}
$$
where $P(n)$ is the set of ordered pairs $(\si,I)$ such that $\si$ is
a permutation of $\{1,2,\dots,n\}$ and $I\succcurlyeq C(\si)$.  For
$(\si,I)\in P(n)$ with $I=(i_1,i_2,\dots,i_l)$ and $0\le k\le l$, let
$\si'_k$ be the permutation of $\{1,2,\dots,n+1\}$ given by
$$
\si'_k(j)=\cases \si(j),& j\le i_1+\dots+i_k\\
n+1,& j=i_1+\dots+i_k+1 \\
\si(j-1), & j>i_1+\dots+i_k+1.\endcases
$$
Also, for $0\le k\le l$ let $I'_k=(i_1,\dots,i_k,1,i_{k+1},\dots,i_l)$,
and for $1\le k\le l$ let
$I''_k=(i_1,\dots,i_{k-1},i_k+1,i_{k+1},\dots,i_l)$; evidently $(\si'_k,I'_k),
(\si'_k,I''_k)\in P(n+1)$ for all $k$.  By iterated application of (6) we have
$$\multline
I[a_{\si(1)}a_{\si(2)}\cdots a_{\si(n)}]*a_{n+1}=
q^{\sum_{i=1}^n |a_i||a_{n+1}|}
a_{n+1}I[a_{\si(1)}\cdots a_{\si(n)}]+\\
\sum_{k=1}^l q^{\sum_{j=i_1+\dots+i_k+1}^n |a_{\si(j)}||a_{n+1}|}
\left(I'_k[a_{\si'_k(1)}\cdots a_{\si'_k(n+1)}]+
I''_k[a_{\si'_k(1)}\cdots a_{\si'_k(n+1)}]\right) .
\endmultline$$
Hence $\Phi_q(a_1\cdots a_{n+1})$ is the sum over $(\si,I)\in P(n)$ of
$$\multline
q^{\sum_{(i,j)\in\iota(\si'_0)}|a_i||a_j|}I'_0[a_{\si'_0(1)}\cdots
a_{\si'_0(n+1)}]+\\
\sum_{k=1}^l q^{\sum_{(i,j)\in\iota(\si'_k)}|a_i||a_j|}\left(
I'_k[a_{\si'_k(1)}\cdots a_{\si'_k(n+1)}]+
I''_k[a_{\si'_k(1)}\cdots a_{\si'_k(n+1)}]\right)
\endmultline$$
and the conclusion follows by noting that every $(\tau,J)\in P(n+1)$
can be written uniquely as one of $(\si'_k,I'_k)$ or $(\si'_k,I''_k)$
for some $(\si,I)\in P(n)$.
\enddemo
In the case $q=0$, our formula for $\Phi_q(w)$ reduces to
$$
\Phi_0(w)=\sum_{I\in\Cal C(\ell(w))}I[w]=(-1)^{\ell(w)}S(\bar w) ,
$$
and by applying Lemma 2.4 with $f(t)=t/(1-t)$ we see that $\Phi_0$ has
inverse $\Phi_0^{-1}$ given by
$$
\Phi_0^{-1}(w)=\sum_{k=1}^{\ell(w)}\sum_{I\in\Cal C(\ell(w),k)}(-1)^{\ell(w)-k}I[w] .
$$
\par
For any word $w=a_1a_2\cdots a_n$, let $V_w$ be the vector space over $k$
with basis $\{a_{\tau(1)}\cdots a_{\tau(n)}\mid\text{permutations $\tau$}\}$, and let
$\phi_{w,q}:V_w\to V_w$ be $\Phi_q$ followed by projection onto $V_w$.
Then $\phi_{w,q}$ is given by
$$
\phi_{w,q}(a_{\tau(1)}\cdots a_{\tau(n)})=\sum_{\text{permutations $\si$}}
q^{\sum_{(i,j)\in\iota(\si)}|a_{\tau(i)}||a_{\tau(j)}|}a_{\si\tau(1)}
\cdots a_{\si\tau(n)},
$$
and we have the following result.
\proclaim{Lemma 5.3} The linear map $\phi_{w,q}$ as defined above has
determinant
$$
\prod_{k=2}^n \prod\Sb\text{$k$-sets}\\ S\subset\{1,\dots,n\}\endSb
\left(1-q^{2\sum_{i,j\in S}|a_i||a_j|}\right)^{(n-k+1)!(k-2)!} .
$$
\endproclaim
\demo{Proof} Following \cite{\DKKT}, we use Varchenko's theorem \cite{\V}
on determinants of bilinear forms on hyperplane arrangements.
To apply the result of \cite{\V}, we consider the set of hyperplanes
in $\bold R^n$ given by $\Cal H_{ij}=\{(x_1,\dots, x_n)\mid x_i=x_j\}$.
To the hyperplane $\Cal H_{ij}$ we assign the weight $\wt\Cal H_{ij}=
q^{|a_i||a_j|}$.  The edges (nontrivial intersections) of this arrangement
are indexed by subsets $S\subset\{1,2,\dots,n\}$ with two or more
elements:  the edge $E_S$ corresponding to the set $S$ is
$$
\bigcap\{\Cal H_{ij}\mid i,j\in S\}=\{(x_1,\dots,x_n)\mid x_i=x_j\ \text
{for all $i,j\in S$}\} .
$$
The edge $E_S$ has weight
$$
\wt E_{S}=\prod_{i,j\in S} \wt\Cal H_{ij}=q^{\sum_{i,j\in S}|a_i||a_j|}.
$$
The domains (connected components) for this hyperplane arrangement are
indexed by permutations:  $C_{\si}=\{(x_{\si(1)},\dots,x_{\si(n)})\mid
x_1<x_2<\cdots<x_n\}$.  Then the quadratic form $B$ on this
arrangement given by
$$
B(C_{\si},C_{\tau})=\prod\Sb\text{hyperplanes $\Cal H_{ij}$}\\
\text{separating $C_{\si}$ and $C_{\tau}$}\endSb \wt \Cal H_{ij}
=\prod_{(i,j)\in \iota(\si\tau^{-1})}q^{|a_{\tau(i)}||a_{\tau(j)}|} 
$$
has the same matrix as $\phi_{w,q}$.  Hence, by Theorem 1.1 of \cite{\V}
we have
$$
\det \phi_{w,q}=\prod_{\text{edges $E$}}(1-\wt(E)^2)^{n(E)p(E)} ,
$$
where the product is over the edges of the hyperplane arrangement,
and $n(E)$ and $p(E)$ are numbers defined in \S2 of \cite{\V}.
It is easy to see from the definitions that $n(E_S)=(n-k+1)!$ and
$p(E_S)=(k-2)!$ for the edge $E_S$ corresponding to a $k$-set
$S\subset \{1,\dots,n\}$, so the conclusion follows.
\enddemo
\proclaim{Theorem 5.4} Any $q\in k$ that is not a root of unity is
generic (i.e., $\Phi_q$ is an isomorphism when $q$ is not a root of
unity).
\endproclaim
\demo{Proof} Suppose $q$ is not a root of unity.  We shall show that
$\Phi_q^{-1}(w)$ exists for any word $w$ by induction on $\ell(w)$.
Using Lemma 5.2 and the induction hypothesis, to find
$\Phi_q^{-1}(a_1\cdots a_n)$ it suffices to find an element $u$ such
that $$
\Phi_q(u)=a_1a_2\cdots a_n + \text{terms of length $<n$}.
$$
But we can do this by taking $u=\phi_{w,q}^{-1}(w)$, and $\phi_{w,q}$
is invertible by Lemma 5.3.
\enddemo
If $q$ is generic, we can define a comultiplication $\De_q$ on $\frak A$
by requiring that all letters be $\De_q$-primitives and that 
$\De_q$ be a $*_q$-homomorphism, i.e. that
$\De_q(a)=a\otimes 1+1\otimes a$ for all $a\in A$ and 
$\De_q(u*_qv)=\De_q(u)*_q\De_q(v)$ for all $u,v\in\frak A$.
This makes $(\frak A,*_q,\De_q)$ a Hopf algebra.  In fact, as we see
in the next result, it is isomorphic to the concatenation Hopf algebra
$(\frak A,\conc,\de)$, where
$$
\de(w)=\sum_{\text{words $u,v$ of $\frak A$}}(u\ \sh\ v,w^*)u\otimes v .
$$
\proclaim{Theorem 5.5} For generic $q$, $\Phi_q$ is a Hopf algebra
isomorphism from $(\frak A,\conc,\de)$ to $(\frak A,*_q,\De_q)$.
\endproclaim
\demo{Proof} Since $q$ is generic, $\Phi_q$ is an algebra isomorphism.
It suffices to show that $(\Phi_q\otimes\Phi_q)\circ\de=\De_q\circ\Phi_q$
on a set of generators:  but this follows because they agree on the
primitives (elements of $A$), which generate $\frak A$ under conc.
\enddemo
\par
In the next result we record a formula for $\De_q(ab)$
when $q$ is generic.  This may be compared with the corresponding formula in
Example 5.2 of \cite{\DKKT}.
\proclaim{Proposition 5.6} Let $a,b,c\in A$.  For $q$ generic,
$$
\De_q(ab)=ab\otimes 1 + 1\otimes ab +\frac{1}{1+q^{|a||b|}}(a\otimes b
+b\otimes a) .
$$
\endproclaim
\demo{Proof} Apply $\De_q$ to the equation
$$
ab=(1-q^{2|a||b|})^{-1}(a*_qb-q^{|a||b|}b*_qa)-(1-q^{|a||b|})^{-1}[a,b] .
$$
\enddemo
A formula for $\De_q(abc)$ can be derived by applying $\De_q$ to
$$
abc=\left(\phi_{abc,q}^{-1}\right)_{\text{id},\text{id}} a*_qb*_qc +
\left(\phi_{abc,q}^{-1}\right)_{\text{id},(12)} b*_qa*_qc +\dots+
\text{terms of length }\le 2 ,
$$
but it is too complicated to give here (it contains twenty terms).
\par
For the cases $q=1$ and $q$ not a root of unity, we have defined a
Hopf algebra $(\frak A,*_q,\De_q)$ with all elements of $A$
primitive.  It would be of interest to define such a Hopf algebra
structure for all $q$.
\subheading {6. Examples}  As we have already remarked, if $[a,b]=0$ for
all generators $a,b\in A$ then $(\frak A,*)=(\frak A,\sh)$ is the
shuffle algebra
as described in Chapter 1 of \cite{\R} (Note, however, that the grading may
be different).  The $q$-shuffle product $\odot_q$ as defined in
\cite{\DKKT, \S4} is the operation
$*_q=\sh_q$ in this case.  This algebra may also be obtained as a
special case of the constructions of Green \cite{\Gr} and Rosso \cite{\Ro}
involving quantum groups.
To identify Green's ``quantized shuffle algebra''
with our construction, take the ``datum''
to be our generating set $A$, with bilinear form $a\cdot b=|a||b|$ for
$a,b\in A$; then Green's algebra $G(k,q,A,\cdot)$ \cite{\Gr, p. 284} is
our $(\frak A,\sh_q)$, except that Green's algebra is $\bold NA$-graded
rather than $\bold N$-graded.  To obtain our algebra from Rosso's
``exemple fondamental'' \cite{\Ro, \S2.1}, take $V$ to be the vector
space over $k$ generated by $A=\{e_1,e_2,\dots\}$, and let
$q_{ij}=q^{|e_i||e_j|}$.
Here are some other examples.
\par
\it Example 1. \rm Let $A_n=\{z_n\}$ for all $n\ge 1$ and
$[z_i,z_j]=z_{i+j}$.   Then $(\frak A,*)$ is just the algebra $\HH$
as presented in \cite{\Ho}.  As is proved there (Theorem 3.4 ff.),
the map $\phi$ defined by
$$
\phi(z_{i_1}z_{i_2}\cdots z_{i_k})=\sum_{n_1>n_2>\cdots>n_k\ge 1}
t_{n_1}^{i_1} t_{n_2}^{i_2}\cdots t_{n_k}^{i_k}
$$
is an isomorphism of $\HH$ onto the algebra of
quasi-symmetric functions over $k$ (denoted $\Qs$ in
\cite{\MaR}).
For each $n\ge 0$, the monomial quasi-symmetric functions
$M_{(i_1,\dots,i_k)}=\phi(z_{i_k}\cdots z_{i_1})$,
where $(i_1,\dots, i_k)\in\Cal C(n)$, form a vector-space basis for
$\frak A_n$.  For our purposes it is more convenient to
identify $M_{(i_1,\dots,i_k)}$
with $z_{i_1}\cdots z_{i_k}$:  under this identification (which is also
an isomorphism), the notation used above is simplified by the
observation that, for compositions $I\in\Cal C(n,k)$ and $J\in\Cal C(k)$,
$J[M_I]=M_{J\circ I}$.  So, e.g., $S(M_I)=(-1)^{\ell(I)}
\sum_{\bar I\succcurlyeq J}M_J$,
where $\bar I$ is the reverse of $I$.  If we let $\Cal L$ denote the
set of $I$ such that $M_I$ corresponds to a Lyndon word, then Theorem 2.6 says
that $\{M_I\mid I\in\Cal L\}$ generates $\frak A=\Qs$ as an
algebra.  The Hopf algebra structure is that described in \cite{\Eh,\MaR};
the two formulas for its antipode are discussed in \cite{\Eh, \S3}.
\par
For the \it integral \rm Hopf algebra $\Qsz$ of quasi-symmetric 
functions, $\{M_I\mid I\in\Cal C(n)\}$ is a $\bold Z$-module basis for
the elements of degree $n$, but $\{M_I\mid I\in\Cal L\}$ is not an algebra
basis.  Nevertheless, from \cite{\DS,\Sc} $\Qsz$ has 
an algebra basis $\{M_I\mid I\in\Cal L^{\text{mod}}\}$,
where $\Cal L^{\text{mod}}$ is the set of 
``modified Lyndon'' or ``elementary unreachable'' compositions, i.e.
concatenation powers of elements of $\Cal L$ whose parts have greatest
common factor 1.
(There is a bijection of $\Cal L$ onto $\Cal L^{\text{mod}}$ given by sending
$(i_1,\dots,i_l)$ to the
$d$th concatenation power of $(\frac{i_1}{d},\dots,\frac{i_l}{d})$, where
$d$ is the greatest common factor of $i_1,\dots,i_k$.)
Of course $\exp$ cannot be defined
over $\bold Z$ because of denominators.
\par
Another algebra basis for $\Qs$ is given by $\{P_I\mid I\in\Cal L\}$,
where $P_I=\exp(M_I)$.
(These are exactly the elements whose duals $P_I^*=\log^*(M_I^*)$ are
introduced in \cite{\MaR, \S2} as a basis for the dual $\Qsd$;
cf. equations (2.12) of \cite{\MaR} and (5) above.)  Since exp is a
Hopf algebra isomorphism, we have the formulas 
$$
P_I*P_J=\sum_{K\in I \sh J} P_K,\quad
\De(P_K)=\sum_{I\sqcup J=K}P_I\otimes P_J,\quad\text{and}\quad
S(P_I)=(-1)^{\ell(I)} P_{\bar I},
$$
where, for compositions $I$ and $J$,
$I\sh J$ is the multiset of compositions obtained by ``shuffling''
$I$ and $J$ (e.g. $(1,2)\ \sh\ (2)=\{(2,1,2),(1,2,2),(1,2,2)\}$), and
$I\sqcup J$ is the concatenation of $I$ and $J$.
\par
Following Gessel \cite{\Ge}, there is still another basis 
$\{F_I\mid I\in\Cal L\}$ for $\Qs$, where
$F_I = \sum_{J\succcurlyeq I} M_J$.
(Then $M_I=\sum_{J\succcurlyeq I} (-1)^{\ell(J)-\ell(I)}F_J$, and since the
coefficients are integral $\{F_I\mid I\in\Cal L^{\text{mod}}\}$ is
a basis for $\Qsz$).
The expansion of the product $F_I*F_J$ in terms of the $F_K$
can be described using
permutations and their descent compositions; see \cite{\TU} or
\cite{\MaR}.  Dualizing Proposition 3.13 and Corollary 3.16 of \cite{\GKLLRT}
(see below), we have
$$
\De(F_K)=\sum_{I\sqcup J=K}F_I\otimes F_J+\sum_{I\vee J=K}F_I\otimes F_J
\quad\text{and}\quad
S(F_I)=(-1)^{|I|}F_{I^{\sptilde}},
$$
where $I\vee J=(i_1,\dots,i_{k-1},i_k+j_1,j_2,\dots,j_l)$ for nonempty
compositions $I=(i_1,\dots,i_k)$ and $J=(j_1,\dots,j_l)$, and
$I^{\sptilde}$ is the conjugate composition of $I$ (as defined in
\cite{\GKLLRT, \S3.2}).  By dualizing Corollary 4.28 of \cite{\GKLLRT}
we have a formula for $F_I$ in terms of the $P_I$:
$$
F_I=\sum_{|J|=|I|}\phr(I,J)\frac{P_J}{\Pi(J)}.
$$
Here $\Pi(I)$ is the product of the parts of the composition $I$, and
$\phr(I,J)$ is as defined in \cite{\GKLLRT, \S4.9}:  for compositions
$I$ and $J=(j_1,\dots,j_s)$ of the same weight, let
$I=I_1\bullet I_2\bullet\cdots\bullet I_s$
be the unique decomposition of $I$ such that $|I_i|=j_i$ for $1\le i\le s$
and each symbol $\bullet$ is either $\sqcup$ or $\vee$; then
$$
\phr(I,J)=\prod_{i=1}^s\frac{(-1)^{\ell(I_i)-1}}{\binom{|I_i|-1}{\ell(I_i)-1}}.
$$
\par
The dual Hopf algebra $\Qsd$ is described in \cite{\MaR, \S2}; it is also
the algebra $\bold{Sym}$ of noncommutative symmetric functions as
defined in \cite{\GKLLRT}. (The coproduct $\de'$ of \S4 corresponds
to the coproduct denoted $\ga$ in \cite{\MaR} and \cite{\GKLLRT}.)
The $M_I$ are dual to the ``products of complete
homogeneous symmetric functions'' $S^I$ (i.e., $(M_I,S^J)=\de_{IJ}$),
while the ``products of power sums of the second kind'' $\Phi^I$
are dual to the elements $P_I/\Pi(I)$ (see \cite{\GKLLRT, \S3} for
definitions).  The $F_I$ are dual to the
``ribbon Schur functions'' $R_I$ \cite{\GKLLRT, \S6}.
\par
The deformation $(\frak A,*_q)$ is the algebra of quantum quasi-symmetric
functions as defined in \cite{\TU}.  The multiplication rule for
``quantum quasi-monomial functions'' as given in \cite{\TU, p. 7345} can
be recognized as (6).
\par
\it Example 2. \rm  For a fixed positive integer $r$, let
$A_n=\{z_{n,i}\mid\ 0\le i\le r-1\}$
and $[z_{n,i},z_{m,j}]=z_{n+m,i+j}$, where the second subscript is
to be understood mod $r$.  By Theorem 2.6, $(\frak A,*)$ is the
polynomial algebra on the Lyndon words in the $z_{i,j}$; by Proposition 2.7,
the number of Lyndon words in $\frak A_n$ is
$$
L_n=\frac1n\sum_{d|n}\mu\left(\frac{n}{d}\right)(r+1)^d
$$
for $n\ge 2$ (and $L_1=r$).  In this case, we call the Hopf algebra
$(\frak A,*,\De)$ the Euler algebra $\EEr$.  Of course $\EE_1$ is the
preceding example (We write $z_i$ for $z_{i,0}$ if $r=1$);
in general there is a homomorphism $\pi_r:\EEr\to\EE_1$
given by $\pi_r(z_{i,j})=z_i$.  
The map $\phi:\EEr\to\calgi$ with
$$
\phi(z_{i_1,j_1}z_{i_2,j_2}\cdots z_{i_k,j_k})=
\sum_{n_1>n_2>\dots>n_k\ge 1}e^{\frac{2\pi i}{r}(n_1j_1+\dots+n_kj_k)}
t_{n_1}^{i_1}\cdots t_{n_k}^{i_k} 
\tag7$$
is an isomorphism of $\EEr$ onto a subring of $\calgi$ (for proof see
\S7 below.)  If we define $\psi_r:\calgi\to\calgi$ by
$$
\psi_r(t_i)=\cases 0,& r \nmid i\\
t_{j}, &i=rj\endcases
$$
(Note $\psi_r$ takes $\Qs\subset\calgi$ isomorphically onto itself!),
then $\psi_r\circ\phi=\phi\circ\pi_r$.  The sets $L$ of Lyndon words in
the $z_{i,j}$ and $\{\exp(w)\mid w\in L\}$ are both algebra bases
for $\EEr$, corresponding to the bases $\{M_I\mid I\in\Cal L\}$ and
$\{P_I\mid I\in\Cal L\}$, respectively, of Example 1.  If we set
$\hat w=\sum_{v\in\Cal P(w)} v$,
where $\Cal P(w)$ is as defined at the end of \S4, then there is a
a basis $\{\hat w\mid w\in L\}$ corresponding to $\{F_I\mid I\in\Cal L\}$.
Note, however,
that while $\pi_r$ maps words to the $M_I$ and exponentials of
words to the $P_I$ (exp commutes with $\pi_r$), in general
$\pi_r(\hat w)$ is not of the form $F_I$.
\par
The dual $\EEr^*$ of the Euler algebra is the concatenation algebra on
elements $z_{i,j}^*$, with coproduct $\de'$ as described in \S4.  The
transpose of $\pi_r$ is the homomorphism $\pi_r^*:\EE_1^*\to\EEr^*$ with
$\pi_r^*(z_i^*)=\sum_{j=1}^{r-1}z_{i,j}^*$.
\par
The motivation for the Euler algebra $\EEr$ comes from numerical series
of the form
$$
\sum_{n_1>n_2>\dots>n_k\ge 1}\frac{\ep_1^{n_1}\ep_2^{n_2}\cdots\ep_k^{n_k}}
{n_1^{i_1} n_2^{i_2}\cdots n_k^{i_k}} ,
\tag8$$
where the $\ep_i$ are $r$th roots of unity and $i_1,i_2,\dots, i_k$
are positive integers (with $\ep_1i_1\ne 1$, for convergence).
In fact (8) is
$\lim_{n\to\infty}\phi_n(z_{i_1,j_1}\cdots z_{i_k,j_k})(1,2,\dots,\frac1{n})$,
where $\phi_n$ is as defined in \S7 and the $j_s$ are chosen appropriately,
so the algebra of such series can be seen as a homomorphic image of
(a subalgebra of) $\EEr$.
These series are called ``Euler sums'' in \cite{\BBB,\Br} and ``values of multiple
polylogarithms at roots of unity'' in \cite{\Go}; in the case $r=1$
the corresponding series are known as ``multiple harmonic series'' \cite{\Ho}
or ``multiple zeta values'' \cite{\Z}.
\par
\it Example 3. \rm Fix a positive integer $m$ and let $A_n=\{z_n\}$ for
$n\le m$ and $A_n=\emptyset$ for $n>m$.  Define
$$
[z_i,z_j]=\cases z_{i+j}&\text{if $i+j\le m$,}\\0&\text{otherwise.}
\endcases
$$
Then $(\frak A,*)$ is the algebra of quasi-symmetric functions on variables
$t_1,t_2,\dots$ subject to the relations $t_i^{m+1}=0$ for all $i$.
\par
\it Example 4. \rm  Let $P(n)$ be the set of partitions of $n$ and let
$A_n=\{z_{\la}\mid\la\in P(n)\}$.  Define $[z_{\la},z_{\mu}]=z_{\la\cup\mu}$,
where $\la\cup\mu$ is the union $\la$ and $\mu$ as multisets.  Then $(\frak A,*)$
can be thought of as the algebra of quasi-symmetric functions in the
variables $t_{i,j}$, where $|t_{i,j}|=j$, in the following sense.  For a
partition $\la=(n_1,\dots,n_l)$, let $t_i^{\la}=t_{i,n_1}\cdots t_{i,n_l}$.
Then any monomial in the $t_{i,j}$ can be written in the form $t_{i_1}^{\la_1}\cdots
t_{i_k}^{\la_k}$, and we call a formal power series quasi-symmetric when the
coefficients of any two monomials $t_{i_1}^{\la_1}\cdots t_{i_k}^{\la_k}$
and $t_{j_1}^{\la_1}\cdots t_{j_k}^{\la_k}$ with $i_1<\cdots<i_k$ and
$j_1<\cdots<j_k$ are the same.
\subheading{7. The Euler algebra as power series} Fix a positive integer
$r$, and let $\EEr$ and $\pi_r:\EEr\to\EE_1$ be as in Example 2.
We shall show $\EEr$ can be imbedded in the formal power series ring
$\bold C[[t_1,t_2,\dots]]$.
For positive integers $n$, define a map $\phi_n:\EEr\to
\calg$ as follows.  Let $\phi_n$ send $1\in\EEr$ to $1\in\calg$ and any
nonempty
word $w=z_{i_1,j_1}z_{i_2,j_2}\dots z_{i_k,j_k}$ to the the polynomial
$$
\sum_{n\ge n_1>n_2>\cdots>n_k\ge 1}\om^{j_1n_1+j_2n_2+\dots+j_kn_k}
t_{n_1}^{i_1}t_{n_2}^{i_2}\cdots t_{n_k}^{i_k},
$$
where $\om=e^{\frac{2\pi i}{r}}$
(If $k > n$, the sum is empty and we assign it the value 0).
Extend $\phi_n$ to $\EEr$ by linearity.  If we make $\calg$ a graded
algebra by setting $|t_i|=1$, then $\phi_n$ preserves the grading.
Also, it is immediate from the definition that
$$
\phi_n(z_{p,i}w)=\sum_{n\ge m > 1}\om^{im}t_m^p\phi_{m-1}(w) 
\tag9$$
for any nonempty word $w$.
\proclaim{Theorem 7.1} For any $n$, $\phi_n:\EEr\to\calg$ is a
homomorphism of graded $k$-algebras.
\endproclaim
\demo{Proof} It suffices to show $\phi_n(w_1*w_2)=\phi_n(w_1)\phi_n(w_2)$
for words $w_1$, $w_2$.  This can be done by induction on
$\ell(w_1)+\ell(w_2)$, following the argument of \cite{\Ho, Theorem 3.2}
(and using equation (9) above in place of equation (*) of \cite{\Ho}).
\enddemo
\proclaim{Lemma 7.2} For $0\le j_1,j_2,\dots,j_m\le r-1$, let
$\cmult\in\bold Q$ be such that
$$
\sum_{j_1=0}^{r-1}\sum_{j_2=0}^{r-1}\cdots\sum_{j_m=0}^{r-1}
\cmult\om^{n_1 j_1+n_2j_2+\dots+n_mj_m}=0 
$$
for all $mr\ge n_1 > n_2 >\cdots> n_m\ge 1$, where $\om=e^{\frac{2\pi i}{r}}$
as above.  Then all the $\cmult$ are zero.
\endproclaim
\demo{Proof} We use induction on $m$.  For $m=1$ the hypothesis is
$$
\sum_{j=1}^{r-1} c_j\om^{nj} = 0 \quad\text{for all}\quad 1\le n\le r,
$$
which is evidently equivalent to having the equality for $0\le n\le r-1$.
But then the conclusion follows from the nonsingularity of the
Vandermonde determinant of the quantities $1,\om, \om^2,\dots, \om^{r-1}$.
\par
Now let $m>1$, and fix $(m-1)r\ge n_2>n_3>\dots>n_m\ge 1$.  Then the
hypothesis says
$$
\sum_{j_1=0}^{r-1}\left(\sum_{j_2=0}^{r-1}\cdots\sum_{j_m=0}^{r-1}\cmult
\om^{n_2j_2+\dots+n_mj_m}\right)\om^{n_1j_1} = 0
\quad\text{for $(m-1)r< n_1\le mr$} .
$$
This is evidently equivalent to having the equality hold for all
$1\le n_1\le r$:  but then we are in the situation of the preceding
paragraph and so
$$
\sum_{j_2=0}^{r-1}\cdots\sum_{j_m=0}^{r-1}\cmult
\om^{n_2j_2+\dots+n_mj_m}=0,
$$
from which the conclusion follows by the induction hypothesis.
\enddemo
\proclaim{Theorem 7.3} The homomorphism $\phi_{nr}$ is injective through
degree $n$.
\endproclaim
\demo{Proof} Suppose $u\in\ker\phi_{nr}$ has degree $\le n$.  Without
loss of generality we can assume $u$ is homogeneous, and in fact that
$\pi_r(u)$ is a multiple of $z_{i_1} z_{i_2}\cdots z_{i_m}$ for $m\le n$.
Then $u$ has the form
$$
u=\sum_{j_1=0}^{r-1}\sum_{j_2=0}^{r-1}\cdots\sum_{j_m=0}^{r-1}
\cmult z_{i_1,j_1} z_{i_2,j_2} \cdots z_{i_m,j_m} ,
$$
and $u\in\ker\phi_{nr}$ implies that
$$
\sum_{j_1=0}^{r-1}\sum_{j_2=0}^{r-1}\cdots\sum_{j_m=0}^{r-1}
\cmult\om^{n_1j_1 + n_2j_2 +\dots+ n_mj_m} = 0
$$
for all $nr\ge n_1>n_2>\cdots > n_m\ge 1$.  But then $u=0$ by the lemma.
\enddemo
\par
For $m\ge n$, there is a restriction map $\rho_{m,n}:\bold C[t_1,\dots,t_m]
\to\bold C[t_1,\dots,t_n]$ sending $t_i$ to $t_i$ for $1\le i\le n$ and
$t_i$ to zero for $i>n$.  Let $\frak P$ be the inverse limit of the
$\bold C[t_1,\dots,t_n]$ with respect to these maps (in the category
of graded algebras); $\frak P$ is a subring of $\bold C[[t_1,t_2,\dots]]$.
The $\phi_n$ define a homomorphism $\phi:\EEr\to\frak P$, and 
the following result is evident.
\proclaim{Theorem 7.4} The homomorphism $\phi$ is injective, and satisfies
equation (7).
\endproclaim

\enddocument